\documentclass[11pt]{article}
\usepackage[utf8]{inputenc}
\usepackage[T1]{fontenc}
\usepackage{lmodern}  
\usepackage{microtype}  
\usepackage{ulem} 

\usepackage[english]{babel}
\usepackage{csquotes}  

\usepackage{amsmath, amssymb, amsthm}
\usepackage{mathtools}  
\usepackage{mathrsfs}   
\usepackage{bbm}        
\usepackage{bm}         
\usepackage{xcolor}  
\newtheorem{theorem}{Theorem}[section]
\newtheorem{lemma}[theorem]{Lemma}
\newtheorem{corollary}[theorem]{Corollary}

\theoremstyle{definition}
\newtheorem{definition}[theorem]{Definition}
\newtheorem{example}[theorem]{Example}
\theoremstyle{remark}
\newtheorem{remark}[theorem]{Remark}

\usepackage{comment}

\usepackage[a4paper, margin=1in]{geometry}
\usepackage{setspace}
\usepackage[colorlinks=true, allcolors=blue]{hyperref}

\newcommand{\eps}{\varepsilon}

\newcommand{\debi}[1]{\textbf{\textit{#1}}} 
\newcommand{\dbi}[1]{\textbf{\textit{#1}}} 
\newcommand{\full}{full}  

\newcommand{\tile}{tile}  
\newcommand{\tiles}{tiles}
\newcommand{\Tile}{Tile}

\newcommand{\meld}{meld}    
\newcommand{\melds}{melds}  

\title{Bootstrap Percolation, Indecomposable Permutations, and the $n$-Kings problem}
\author{Mark Huibregtse, Cristobal Lemus-Vidales, and David Vella}
\date{\today}

\begin{document}

\maketitle
 
\begin{abstract}
    We study the process of bootstrap percolation on $n \times n$ permutation matrices, inspired by the work of Shapiro and Stephens \cite{SS}. In this percolation model, cells mutate (from 0 to 1) if at least two of their cardinal neighbors contain a 1, and thereafter remain unchanged; the process continues until no further mutations are possible.  After carefully analyzing this process, we consider how it interacts with the notion of (in)decomposable permutations.  We prove that the number of indecomposable permutations whose matrices ``fill up'' to contain all 1's (or are ``full'') is half of the total number of full permutations.  This leads to a new proof of a key result in \cite{SS}, that the number of full $n \times n$ permutations is the $(n-1)^{\text{st}}$ large Schr\"{o}der number.  Finally, after rigorously justifying a heuristic argument in [5], we find a new formula for the number of $n \times n$ ``no growth'' permutation matrices, and hence a new solution to the well-known $n$-kings problem. 
\end{abstract}
\section{Introduction}

\subsection{Motivation}

This article was inspired by a close reading of Shapiro and Stephens' paper \cite{SS}, which studied the combinatorics of bootstrap percolation, a growth model in which the value of a cell in a 0-1 matrix can change from 0 to 1 whenever at least two of its neighbors (to the north, south, east, and/or west) have the value 1.  After a cell has ``mutated'' from 0 to 1, its value remains constant thereafter.  The process begins with an $n\times n$ permutation matrix, and proceeds in a series of steps consisting of the random selection of a ``mutable'' cell (if any) and its mutation; the process terminates when there are no more mutable cells.

Two possible outcomes of the process, illustrated by the following examples, are of particular interest:

\begin{example}\label{exmp:FullVsNoGrowth}
The matrix ``fills up'' (each arrow represents two steps):
\[
    \begin{array}{ccccccc}
       \left( \begin{array}{ccc}
                   0 & 0 & 1\\
                   1 & 0 & 0\\
                   0 & 1 & 0
              \end{array}\right)
        & \mapsto &
        \left( \begin{array}{ccc}
                   0 & 0 & 1\\
                   1 & 1 & 0\\
                   1 & 1 & 0
              \end{array}\right)
        & \mapsto & 
        \left( \begin{array}{ccc}
                   1 & 1 & 1\\
                   1 & 1 & 0\\
                   1 & 1 & 0
              \end{array}\right)
        & \mapsto & 
        \left( \begin{array}{ccc}
                   1 & 1 & 1\\
                   1 & 1 & 1\\
                   1 & 1 & 1
              \end{array}\right)        
    \end{array}.
\]
The matrix is ``no-growth'':\\\
\[
  \left( \begin{array}{cccc}
                   0 & 1 & 0 & 0\\
                   0 & 0 & 0 & 1\\
                   1 & 0 & 0 & 0\\
                   0 & 0 & 1 & 0
              \end{array}\right).
\]
\end{example}
Indeed, Shapiro and Stephens set out to determine the number and asymptotic density of the $n\times n$ permutation matrices that 1) fill up, and 2) are no-growth.  In the former case, using generating function techniques, they showed \cite[Theorem 1, p. 276]{SS} that:
\begin{equation}\label{eqn:ExmpFullVsNoGrowth}
    \begin{array}{l}
    \text{\textit{The number of $n\times n$ permutation matrices that fill up is $S_{n-1}$,}}\\
    \text{\textit{the $(n-1)^{\text{st}}$ (big) Schr\"{o}der number.}}
    \end{array}
\end{equation}
From this it follows that the asymptotic density $(S_{n-1})/n!\rightarrow 0$ as $n\rightarrow \infty$.

The Schr\"{o}der numbers (also called large or big Schr\"{o}der numbers): 1, 2, 6, 22, 90, 394, 1806, 8558, 41586,\dots count, among other objects, the number of subdiagonal lattice paths from $(0, 0)$ to $(n, n)$ consisting of steps East $(1, 0)$, North $(0, 1)$ and Northeast $(1, 1)$. See \cite[A006318]{oeis} for more information about this sequence. 

In the second case, they obtained a fascinating answer in terms of a composition of generating functions.  For $n\geq 0$, let $a_n$ denote the number of $n\times n$ no-growth permutation matrices (with $a_0 = 1$), and let $A(x)=\sum_{n=0}^{\infty}a_n x^n$ be the generating function of this sequence.  Let $R(x)$ denote the generating function for the Schr\"{o}der numbers; hence $x R(x)$ is the generating function for the sequence $(S_{n-1})_{n=1}^{\infty}$.  Then Shapiro and Stephens state that 
\begin{equation}\label{eqn:SSCompFormula}
    A(x R(x)) = \sum_{n=0}^{\infty}n! x^n = \eps(x),
\end{equation}
justifying the result with a very brief intuitive summary.  They then note that the functional inverse of $x R(x)$ is $x((1-x)/(1+x))$, yielding
\begin{equation}\label{eqn:keyCompForAn}
    A(x) = \eps(x((1-x)/(1+x))),
\end{equation}
which, as they said, ``at least in some sense,'' provides a computation of $a_n$ and permits the asymptotic density of the no-growth permutation matrices to be deduced.
It was this formula that particularly sparked our interest in \cite{SS}, since a result of D. Vella \cite{V} on composite generating functions would permit the extraction of an explicit formula for the no-growth numbers: 1, 1, 0, 0, 2, 14, 90, 646, 5242,\dots The sequence $a_n$ appears on the  OEIS as \cite[A002464]{oeis}.

Before writing up this new explicit formula (it appears below in the concluding Section \ref{sec:NKingsProb}), we wanted to make sure we thoroughly understood the proof of (\ref{eqn:SSCompFormula}).  This spurred us to undertake a close reading of Shapiro and Stephens's paper. 

The next section discusses the new results we obtained and lays out the overall structure of this paper.  Before beginning this survey, we note a few conventions that are used throughout the paper:  
\begin{itemize}
    \item The set $\{1, 2, \dots, n\}$ of positive integers less than or equal to $n$ is denoted $[n]$.
    \item A permutation $\pi=a_1 a_2 \dots a_n$ of $[n]$ is written in \textbf{\textit{single-line notation}}, that is, given $\pi:[n]\rightarrow [n]$, we have $a_i=\pi(i)$ for $1\leq i \leq n$.
\end{itemize} 

\subsection{Results and Organization}

In Section \ref{sec:BootstrapPerc}, we discuss the process of bootstrap percolation on permutation matrices in detail. We establish some elementary but useful results, including  Theorem \ref{thm:InvarProp}, which states that the final result of percolation does not depend on the order in which mutable cells are changed from $0$ to $1$, and Theorem \ref{thm:finalConfig} and its corollary, which precisely describe the possible final configurations of percolation on a permutation matrix.  

We then turn in \S \ref{subsec:left-merging} to a study of ``\tile\ merging'' algorithms for percolation.  The first part of Example \ref{exmp:FullVsNoGrowth} gives the idea: We first merge two $1\times 1$ \tiles\ (the $(2,1) \text{ and } (3,2) $ cells) to obtain a $2\times 2$ \tile, which is then merged with the $1\times 1$ \tile\ (at (1,3)) to yield the final $3\times 3$ \tile.  Looking at each \tile\ merge as a single step, one achieves greater efficiency of execution and conceptual clarity.  

\Tile\ merging is not mentioned explicitly in \cite{SS}, but is implicit in the corresponding idea of ``bracketing'' of permutations, discussed in \S \ref{subsec:Bracketing}.  The basic idea of bracketing is to encode the sub-\tiles\ that are created during a \tile\ merging sequence with subsequences of the permutation: each sub-\tile\ $T$ corresponds to the initial $1$'s in the permutation matrix that belong to $T$.  In the preceding example, the initial permutation is $\pi=213$, and the result of what we call ``left-bracketing'' is $([21]3)$. This records that the first \tile\ merge joined the $1\times 1$ \tiles\ representing the 1's corresponding to $a_1 = 2 \text{ and } a_2 = 1$, and the second joined the new $2\times 2$ \tile\ with the $1\times 1$ \tile\ representing the 1 corresponding to $a_3=3$. 

If a permutation matrix associated to a permutation $\pi$ fills up under bootstrap percolation, or is ``\full,'' as we shall say, then \tile\ merging/bracketing must lead to a final $n\times n$ \tile\ or single ``top-level bracketing'' having one of the forms $(m_1, m_2)$ or $[m_1,m_2]$; this is discussed in \S \ref{subsec:TopLevelBracketings}.  
This bracketing carries valuable information regarding the underlying permutation, as will be seen. 

It turns out to be important to be precise in describing an algorithm for \tile\ merging or (equivalently) bracketing.  The semi-precise bracketing process discussed in \cite[p.\ 277]{SS} is in some cases compatible with multiple interpretations that can lead to errors if the wrong interpretation is used; this is explained in \S \ref{subseq:BinaryTree}, where we also state and prove a fully valid version of the potentially erroneous assertion in \cite{SS} that could result from such a misinterpretation.

In \S \ref{sec:IndecPerms}, we begin consideration of an issue not discussed in \cite{SS}, namely, the interplay of the concepts of (in)decomposable permutation and bootstrap percolation.  \S \ref{subseq:IndecPermDef} recalls the definitions, and \S \ref{subsec:DefOfIndecComps} discusses the indecomposable components of a given permutation and how to find them.  This discussion is extended in \S \ref{sec:IndecCompOfFullPermutations}, where we consider (in)decomposable permutations that are \full.  In \S \ref{subsec:InvOnFullPerms}, we show that the operation of inversion on permutations ($a_1 a_2 \dots a_n \mapsto a_n a_{n-1} \dots a_1$) carries \full\ permutations to \full\ permutations; we also show that a \full\ indecomposable (resp.\ decomposable) permutation has top-level bracketing (under any bracketing algorithm) of the form $[m_1, m_2]$ (resp.\ $(m_1, m_2)$).  Consequently, since (as we show) inversion changes the form of the top-level bracketing of a full permutation ($[\ ,\, ] \text{ to } (\ ,\, )$ and vice versa),  in \S \ref{subseq:halfLemmaSec} we obtain a key result (Corollary \ref{cor:HalfLemma}, the ``Half-lemma''): for $n\geq 2$, the number of \full\ indecomposable permutations of $[n]$ is exactly half of the number of \full\ permutations.  

We then consider the indecomposable components of a \full\ permutation in \S \ref{subsec:IndecCompOfFullPermutations}. Theorem \ref{thm:CompsThruBracketing} shows how the rightmost indecomposable component of a \full\ permutation can be computed through (left) bracketing; iteration (Example \ref{exmp:indecompComps}) then leads to a process for extracting a permutation's full indecomposable factorization through (left) bracketing. This section concludes with Theorem \ref{thm:piFullIFFCompsFull}, which establishes that a permutation $\pi$ is \full\ if and only if each of its indecomposable components is \full, a result needed in the sequel.

In \S \ref{sec:CountFullPerms}, we offer a new proof of Shapiro and Stephens' result (\ref{eqn:ExmpFullVsNoGrowth}) that the number $p_n$ of \full\ permutations of $[n]$ equals the $(n-1)^{\text{st}}$ (big) Schr\"{o}der number.  Instead of using generating function techniques, our proof shows (using the half-lemma) that the sequence $(p_n)_{n\geq 1}$ satisfies the same recurrence as does the Schr\"{o}der sequence.

Finally, to bring the paper to a close, we offer a detailed proof of equation (\ref{eqn:SSCompFormula}), and then show how it leads to an explicit formula for the number $a_n$ of no-growth permutations of $[n]$.  Since the $n\times n$ no-growth permutation matrices correspond to the arrangements of $n$ non-attacking kings on the $n \times n$ chessboard, we thereby attain an explicit solution to the well-known $n$-kings problem.

\section{Bootstrap Percolation}\label{sec:BootstrapPerc}

This section studies the process of bootstrap percolation on permutation matrices.  It includes a proof that the final configuration attained does not depend on the order in which 
``mutable'' cells are changed from $0 \text{ to } 1$, and a discussion, from our point of view, of the ``bracketing'' process (an efficient means for determining the final configuration) introduced in \cite{SS}.

\subsection{The process of bootstrap percolation}

Let $M$ be an $n\times n$ permutation matrix (a matrix whose entries are $0$'s and $1$'s and in which every row and every column contains exactly one $1$).  Bootstrap percolation is an iterative process in which certain ``mutable'' entries (or cells) of $M$ are changed from $0$ to $1$. 
\begin{definition}\label{def:Mutable}
A cell of $M$ is \debi{mutable} if it contains a $0$ and at least two of its neighbors (to the north, south, east, and/or west) contain a $1$; when the cell is changed, we say it has \debi{mutated}, and is therefore no longer mutable.
\end{definition}
Each step in bootstrap percolation consists of choosing a mutable cell (if any) in the current configuration of the matrix and mutating it.  If several mutable cells are available at any stage, the choice of which to select and mutate is arbitrary.  Percolation continues until there are no remaining mutable cells, at which point it terminates in a final configuration.  The following example illustrates the process:
\begin{example}\label{exmp:BootstrapExamples}
\[
    \begin{array}{ccccc}
      \left(
        \begin{array}{ccccc}
        0 & 0 & 0 & 1 & 0\\
        0 & 1 & 0 & 0 & 0\\
        1 & 0 & 0 & 0 & 0\\
        0 & 0 & 0 & 0 & 1\\
        0 & 0 & 1 & 0 & 0
        \end{array}
      \right)
      &
        \mapsto
      &
      \left(
        \begin{array}{ccccc}
        0 & 0 & 0 & 1 & 0\\
        0 & 1 & 0 & 0 & 0\\
        1 & 1 & 0 & 0 & 0\\
        0 & 0 & 0 & 0 & 1\\
        0 & 0 & 1 & 0 & 0
        \end{array}
      \right)
       &
        \mapsto
      &
      \left(
        \begin{array}{ccccc}
        0 & 0 & 0 & 1 & 0\\
        1 & 1 & 0 & 0 & 0\\
        1 & 1 & 0 & 0 & 0\\
        0 & 0 & 0 & 0 & 1\\
        0 & 0 & 1 & 0 & 0
        \end{array}
      \right)
    \end{array}
\]
\end{example}

\begin{remark}\label{rem:permToPermMatrix}
    The preceding example illustrates our convention for associating a permutation to a permutation matrix.  We view the initial permutation matrix as the ``graph'' of the permutation $34152$, that is, 
\[
1\mapsto 3,\ 2\mapsto 4,\ 3\mapsto 1,\ 4\mapsto 5,\ \text{and}\ 5\mapsto 2.
\]
\end{remark}
In the next two sections, we show that the final result attained by bootstrap percolation is independent of the order in which mutable cells are mutated, and we describe the  properties of the final configuration.

\subsection{Invariance of final result of bootstrap percolation}

The purpose of this section is to prove that the final configuration of bootstrap percolation is independent of the order in which the percolation is carried out (i.e., the order in which mutable cells are mutated).

Given an $n\times n$ permutation matrix $M$, there are $n$ cells containing a 1 and $n^2-n$ cells containing a 0. It follows that bootstrap
percolation will require at most $n^2-n$ steps to reach a terminal position. 

\begin{definition}\label{def:PotMutable}
Given a cell $c=c_{i,j}$ in $M$, we say that $c$ is \debi{potentially mutable} if $c$ initially contains a $0$ and there is a sequence of percolation steps that mutates it.  We then define  
\[
    \begin{array}{rcl}
   L(c) & = & \left\{
        \begin{array}{l}
   \text{the minimum number of steps required to mutate $c$, if $c$ is}\\
   \text{potentially mutable, and otherwise}\\
   n^2 \text{ (which exceeds the maximum possible number of steps).} 
       \end{array}\right.
    \end{array}
\]
We say that a percolation
sequence is {\bf  complete} if it ends in a position with no mutable cells, and therefore no further steps are possible.
Finally, we define a collection of subsets of cells in $M$ as follows: \\
\[
        \begin{array}{rcl}
U_0 &  =  & \{c:c \text{ is initially set to $1$} \}, \text{ and for } 1 \leq i \leq n^2 -n,\\ 
  U_i &  =  & \{c:L(c)=i\}
    \end{array}.
\]
\end{definition}

Clearly $\left|U_0\right|=n$; note it is possible for all the $U_i$ for $1 \leq i \leq n^2-n$ to be empty, which occurs if the original permutation matrix $M$ is in a no-growth configuration (meaning no cell in $M$ is mutable, or, equivalently, $U_1= \emptyset$) .

\begin{example}Consider the following $3 \times 3$ permutation matrix representing the permutation $\pi= 213$:
\[
\left(
\begin{array}{ccc}
 0 & 0 & 1 \\
 1 & 0 & 0 \\
 0 & 1 & 0 \\
\end{array}
\right).
\]
In this matrix, one sees by inspection that 
\[
  \begin{array}{c}
  U_1=\left\{c_{\{2,2\}},c_{\{3,1\}}\right\},\ U_2=\left\{c_{\{1,2\}},c_{\{2,3\}}\right\},\ U_3 = \left\{c_{\{1,1\}},c_{\{3,3\}}\right\}, \text{ and }\vspace{.05in}\\ U_4 = \emptyset,\ U_5 = \emptyset,\ \text{and } U_{6} = U_{3^2 -3} = \emptyset.
  \end{array}
\]
\end{example}

We now turn to the invariance proof, for which we need the following Lemma:
\begin{lemma}\label{lem:InvarLem} If a cell $c$ containing a $0$ is mutable before a percolation step is carried out, and that step does not change its value, it remains mutable after the step is carried out.
\end{lemma}

\proof The number of neighbors of $c$ containing a 1 either stays the same or increases. \qed

\begin{theorem}\label{thm:InvarProp} Any complete sequence of percolation steps mutates the same set of cells as any other, which implies that the final configuration attained does not depend on the order in which the percolation steps are carried out.
\end{theorem}
\begin{proof} We will show that any complete sequence of percolation steps has a final configuration in which the set of cells containing
a $1$ is equal to $\displaystyle\bigcup _{ i=0}^{ n^2-n}U_i$.  The desired result is then immediate. 

In case the permutation matrix  $M$ is no-growth, all the $U_i=\emptyset$ except for $U_0$, which is the set of cells containing the $1${'}s
of $M$. In this case, there is only one complete sequence of percolation steps, namely, the sequence of no steps at all, so
the final configuration is indeed given by $\displaystyle\bigcup_{ i=0}^{ n^2-n}U_i$ in this case.

If we are not in the no-growth case, we must have that $U_1\neq \emptyset$, and all the cells $c \in U_1$ are mutable. The first
percolation step must be to mutate one of the cells in $ U_1$. By the Lemma, the remaining cells in $U_1$ remain mutable,
and will continue to remain so (if left unmutated) as percolation continues. Since a complete percolation sequence must continue until there are no more mutable
cells, eventually all the cells in $U_1$ will be mutated. At that point, every element in $U_2$ will either contain a 1 from
a previous step or be mutable, since every cell mutable in one step has been mutated. By the same argument, eventually all of
the cells in $U_2$ will be mutated, at which point any unmutated cells in $U_3$ will be mutable. Continuing in this way, we see that $\displaystyle\bigcup_{ i=0}^{ n^2-n}U_i$ is a subset of the final configuration of $1$'s, but since the maximum number of steps is equal to $n^2-n$, we know that
there cannot be any other $1$'s in the final configuration. This completes the proof.
\end{proof}

\begin{remark}\label{rem:Newman'sLemma}
    Theorem \ref{thm:InvarProp} is an easy consequence of Newman's Lemma \cite{Newman} (see \cite{NewmanWikipediaPage} for a nice discussion); we thank Ira Gessel for this observation.
\end{remark}

\subsection{The final configuration of bootstrap percolation}\label{sec:finalConfig}

In this section we study the final configuration of bootstrap percolation, that is, the pattern of $1$'s in the matrix at the end of a complete percolation sequence. Here and elsewhere the following terminology is useful:

\begin{definition}\label{def:Tile}
    Given a matrix $(c_{i,j})$, a ($k \times l$) \debi{\tile} is a rectangular sub-region of the matrix containing the following cells:
\[
         \left\{c_{i,j}: i_0 \leq i \leq i_0+k-1 \text{ and } j_0 \leq j \leq j_0+l-1\right\}
\]    
If every cell of a tile $T$ contains a $0$ (resp.\ $1$), we call $T$ a \debi{null} (resp.\ \debi{unitary}) \tile. Lastly, if $T$ is a $k\times k$ \tile\, we say $T$ is a \debi{square} \tile\ of \debi{size} $k$.
\end{definition}

The main result of this section can then be summarized as follows: the final configuration of bootstrap percolation on an $n\times n$ permutation matrix consists of $m$ (for some $1\leq m \leq n$) square unitary \tiles\ in a no-growth configuration.
\par
\bigskip
Imagine the percolation process broken down into steps:

\begin{equation}\label{eqn:stringOfMs}
  M_0 \longrightarrow M_1 \longrightarrow M_2 \longrightarrow ... \longrightarrow M_r,
\end{equation}
where $M_0$ is an $n\times n$ permutation matrix, each step involves mutating a single cell, and $M_r$ is the final matrix 
after percolation ends. We have the following
\par
\bigskip
\begin{lemma}\label{lem:MTRowsCols}
    Every row and every column in each $M_t$ contains exactly one unbroken string of nonzero entries. In other words, a row or a column cannot contain any $0$ entries between any two $1$ entries, and there are no rows or columns of all $0$'s. 
\end{lemma} 

\begin{proof}
    We induct on $t$. The result holds when $t=0$ because we are starting with a permutation matrix $M_0$, which has exactly one 
nonzero entry in each row and column. Assume the result holds for some $M_t$, $1\leq t < r$. The only difference between $M_t$ and $M_{t+1}$ is that a 
single cell, say $c = c_{i,j}$, has mutated, meaning that $c$ is mutable in $M_t$. This means at least two entries 
in $M_t$ that are adjacent to $c$ are already $1$. If these two entries were in the same row (or the same column) 
as $c$ (which contains a $0$), with $c$ adjacent to them both, then $c$ must be between them, which would violate the inductive hypothesis ($M_t$ would contain a row or a column with nonzero entries separated by a $0$ entry). It follows that those nonzero entries which are adjacent to $c$ must be in different rows and columns, but since they are both adjacent to 
$c$, they must be diagonally adjacent to each other, one in the same ($i^{\text{th}}$) row as $c$, and the other in the same ($j^{\text{th}}$) column. In this case, the mutated $c$ (now containing a $1$) in $M_{t+1}$ is added to the $i^{\text{th}}$ row, but is adjacent to an existing nonzero 
entry in $M_t$ already in the $i^{\text{th}}$ row, and similarly, it is added to the $j^{\text{th}}$ column but adjacent to a nonzero entry of 
$M_t$ already in the $j^{\text{th}}$ column. Thus each row and column in $M_{t+1}$ is the same as that in $M_t$, except for one row and one column that acquire a new ``$1$'' entry adjacent to a previous ``$1$'' entry in that row or column. Therefore, no 
gaps are produced and the nonzero entries continue to form a single unbroken string in each row and column. So the Lemma holds for $M_{t+1}$, and so for all $M_t$ by induction.
\end{proof}
\par
\bigskip
Applying Lemma \ref{lem:MTRowsCols} to the final matrix $M_r$ leads to the following Theorem: 

\begin{theorem}\label{thm:finalConfig}
The final matrix $M_r$ consists of square unitary \tiles, separated by $0$'s. Furthermore, each row
and column of $M_r$ meets exactly one of these unitary \tiles.
\end{theorem}

\begin{proof}
Consider the final matrix $M_r$. We know the first column has a string of consecutive nonzero entries, and so does 
the second column. Since $M_r$ is the final matrix, the percolation process has terminated. So, the nonzero strings 
in the first two columns either occupy the exact same rows, or else they are separated by enough zeros to prevent 
further percolation. Any other configuration - such as being diagonally adjacent, or partially overlapping - would lead to 
further percolation, a contradiction. In case the strings in the first two columns do match rows exactly, look at 
the unbroken string in the third column. It must occupy the exact same rows as the string in the second (or be 
separated by enough zeros to avoid percolation), etc.  The upshot is, by iteration, for several columns, all the 
unbroken strings are in the exact same rows. The first time one hits a column, say the $(k+1)^{\text{st}}$, where the unbroken 
string is not in those rows, the rest of those rows in $M_r$ are all $0$ (otherwise, we would contradict Lemma \ref{lem:MTRowsCols}.) 
That means that the nonzero entries in the first $k$ columns form an unbroken rectangular unitary \tile.
Again, so as to not violate Lemma \ref{lem:MTRowsCols}, we note that both the rows and columns passing through this \tile\ meet no other nonzero entries outside of this \tile.
\par
\bigskip
At this point, if we haven't filled up the matrix completely, and there is a $(k+1)^{\text{st}}$ column with its unbroken string 
in a different set of rows than the first \tile, we iterate this process, producing a second rectangular unitary \tile, in consecutive columns beginning with the $(k+1)^\text{st}$, with all the nonzero entries in the same set of rows. By iteration, 
this shows $M_r$ consists of rectangular unitary \tiles, separated by $0$'s, and the additional fact that each 
row and column of $M_r$ meets exactly one of these unitary \tiles\ (because of Lemma \ref{lem:MTRowsCols}).
\par
\bigskip
It remains to see these unitary \tiles\ are actually square. Consider one of them and suppose it has size $p \times q$. Looking at the $q$ columns making up this \tile, we 
observe that in the initial matrix $M_0$, each of those columns has a ``$1$'' entry which survives the percolation process 
and ends up in $M_r$, so \textit{must} be within the rows meeting this \tile\ in $M_r$, since the columns meeting this \tile\  
meet no other nonzero entries outside the \tile. However, those original ``$1$'' entries from $M_0$ must all be 
in distinct rows since $M_0$ is a permutation matrix. It follows that $p \geq q$. By making the same argument working 
with rows (or applying the same argument to the transpose matrix), we get that the $p$ ``$1$'' entries from those rows in $M_0$ 
must end up in (distinct) columns inside the \tile, so $q \geq p$.  It follows that $p=q$, so the unitary \tiles\ are, in fact, square, completing the proof of the Theorem.
\end{proof}
\par

Since each column of $M_r$ meets exactly one of the $m$ square unitary \tiles, we can unambiguously order these \tiles\ from left to right, and label them $T_1, T_2, \dots, T_m$.  We denote the size of $T_j$ by $s_j$; it is then a clear consequence of Lemma \ref{lem:MTRowsCols} that $\displaystyle\sum_{j=1}^m s_j=n$ (the dimension of the matrices $M_0, M_1, \dots, M_r$).
\par
\bigskip
To complete the work of this section, we show that, as asserted above, the unitary \tiles\ $T_j$ of $M_r$ are in a no-growth configuration.  The idea is to condense the matrix $M_r$ into an $m\times m$ permutation matrix in which (reading left-to-right) the $j^{\text{th}}$ cell containing a ``$1$''  corresponds to the $j^{\text{th}}$ unitary \tile\ $T_j$ of $M_r$.  Intuitively, one does this by ``viewing the matrix $M_r$ from far away,'' so that the $T_j$ appear to be single cells (or size $1$ square unitary \tiles).  More precisely, one can construct an $m\times m$ matrix by crossing out, for each $1\leq j \leq m$, $s_j-1$ of the rows and $s_j-1$ of the columns of $M_r$ that meet $T_j$.  What remains is an $m\times m$ permutation matrix $CM$ (because each row and each column of $CM$ contains exactly one $1$ and $m-1$ $0$'s) that is a ``condensed'' version of $M_r$: the $j^{\text{th}}$ $1$ in $CM$ (from the left) corresponds to the \tile\ $T_j$ in $M_R$. We now have the following
 
\begin{corollary}
    The $m\times m$ permutation matrix $CM$ is no-growth, that is, it contains no mutable cells.
\end{corollary}

\begin{proof} Suppose to the contrary that there is a mutable cell $c_{i,j}$ in $CM$.  Then $c_{i,j}$ contains a $0$ and has a $1$ as either west or east neighbor and a $1$ as either north or south neighbor, so the two $1$'s are diagonally adjacent.  Consider the case where $c_{i,j-1}=1$ is the west neighbor and $c_{i+1,j}=1$ is the south neighbor of $c_{i,j}$.  Then we know that the west neighbor corresponds to the unitary \tile\ $T_{j-1}$ and the south neighbor to the unitary \tile\ $T_j$.  Since $M_r$ is the final configuration, we know that $T_{j-1}$ and $T_j$ cannot meet diagonally at a corner, otherwise there would be additional percolation steps to perform.  Accordingly, since the rightmost column of $T_{j-1}$ is adjacent to the leftmost column of $T_j$, there must be one or more rows of $M_r$ that separate the southernmost row of $T_{j-1}$ and the northernmost row of $T_j$.  Consider one of these separating rows, which must intersect a unitary \tile\ $T_k$ distinct from $T_{j-1}$ and $T_j$. After crossing out $s_k-1$ of the $s_k$ rows meeting $T_k$ (all of which are separating rows), the remaining row intersecting $T_k$ must condense to a row in $CM$ that vertically separates the cells $c_{i,j-1}$ and $c_{i+1,j}$, which is a contradiction.  Similar contradictions obtain in the other cases, completing the proof.
\end{proof}

A case of particular interest is highlighted in the following

\begin{definition}\label{def:full}
    We say that a permutation $\pi$ of $[n]$ is \debi{full} if the final configuration of bootstrap percolation on the permutation matrix associated to $\pi$ consists of a single unitary \tile\ of size $n$.
\end{definition}

\subsection{\Tile\ merging algorithms for percolation}\label{subsec:left-merging}

Let $\pi = a_1 a_2 \dots a_n$ be a permutation of $[n]$ and let $M_0$ denote the permutation matrix corresponding to $\pi$.  By Theorem \ref{thm:InvarProp}, we know the final configuration $M_r$ of bootstrap percolation on $M_0$ is independent of the sequence in which mutable cells are mutated.  In this section, we discuss algorithms for bootstrap percolation that provide conceptual insight and allow for rapid computation of the final configuration $M_r$.  The basic idea is to replace the mutation of individual cells with the ``merging'' of diagonally-adjacent square unitary \tiles. We will soon see how each \tile-merging algorithm corresponds to a bracketing process on the string $``a_1 a_2 \dots a_n$.''

\begin{example}\label{exmp:Block Merging}
    Here is an example of the merging of two diagonally-adjacent square unitary \tiles\ of sizes $2$ and $1$.  Although the merging requires four mutation steps, we conceive of the merging as a single ``macro'' step, and don't specify the order in which the mutations will occur.  Also note that if two square unitary \tiles\ of sizes $p$ \text{and} $q$ are merged, the result is a square unitary \tile\ of size $(p+q)$.   
\[
    \begin{array}{ccccccccc}
      \left(
        \begin{array}{ccc}
        1 & 1 & 0\\
        1 & 1 & 0\\
        0 & 0 & 1
        \end{array}
      \right)
      &
        \mapsto
      &
      \left(
        \begin{array}{ccc}
        1 & 1 & 0\\
        1 & 1 & 1\\
        0 & 0 & 1
        \end{array}
      \right)
       &
        \mapsto
      &
      \left(
        \begin{array}{ccc}
        1 & 1 & 1\\
        1 & 1 & 1\\
        0 & 0 & 1
        \end{array}
      \right)
     &
       \mapsto
     &  
     \left(
        \begin{array}{ccc}
        1 & 1 & 1\\
        1 & 1 & 1\\
        0 & 1 & 1
        \end{array}
      \right)
      &
        \mapsto
      &
      \left(
        \begin{array}{ccc}
        1 & 1 & 1\\
        1 & 1 & 1\\
        1 & 1 & 1
        \end{array}
      \right)
    \end{array}
\]
\end{example}

In brief, a \debi{\tile\ merging algorithm} proceeds as follows: One begins with the $n \times n$ permutation matrix $M_0$, which contains $n$ square unitary \tiles\ of size $1$. One repeatedly searches the square unitary \tiles\ in the matrix (in some specified order), and merges diagonally-adjacent pairs until no more merges are possible, at which point the algorithm terminates.

On termination, there are two possibilities:
\begin{itemize}
    \item[i.] There is a single square unitary \tile\ of size $n$.  This happens when the original permutation is \full.
    \item[ii.] There are $m$, for $1< m \leq n$, square unitary \tiles\ in $M_r$ that are in a no-growth configuration, since no further merges are possible.
\end{itemize}
It follows that the final configuration has the form described in Section \ref{sec:finalConfig}, that is, it will be a matrix containing $1\leq m \leq n$ square unitary \tiles\ in a no-growth configuration, which implies by definition that a complete sequence of cell mutations has been carried out. Consequently, by Theorem \ref{thm:InvarProp}, any other complete sequence of cell mutations beginning with the original permutation matrix will lead to the same final configuration. 

We are primarily concerned with the algorithms named in the following
\begin{definition}
    The \debi{left merging} algorithm proceeds as follows: The matrix being percolated (starting with $M_0$) is scanned from left to right in search of a diagonally adjacent pair of square unitary \tiles.  If such a pair is found, the \tiles\ are merged and the process is repeated until a traversal finds no diagonally adjacent pair to merge, at which point the process terminates.  It is clear that termination must occur, since the number of square unitary \tiles\ decreases by $1$ each time a merge occurs.  The \debi{right merging} algorithm is essentially the same, except that the matrix being percolated is repeatedly scanned from right to left.
\end{definition}

Here is an example of these algorithms on a $4\times 4$ (\full) permutation matrix; each arrow represents the execution of one traversal and merge operation:

\begin{example}\label{exmp:leftRightBlockMerge}
    Let $\pi = (1324)$, a (full) permutation of $[4]$.  We see that left merging and right merging lead to the same final configuration (as they must), but through different pathways.
    \par
    \textbf{Left merging:}
\[
    \begin{array}{ccccccc}
      \left(
        \begin{array}{cccc}
        0 & 0 & 0 & 1 \\
        0 & 1 & 0 & 0 \\
        0 & 0 & 1 & 0 \\
        1 & 0 & 0 & 0 
        \end{array}
      \right)
      &
        \mapsto
      &
      \left(
        \begin{array}{cccc}
        0 & 0 & 0 & 1 \\
        0 & 1 & 1 & 0 \\
        0 & 1 & 1 & 0 \\
        1 & 0 & 0 & 0         
        \end{array}
      \right)
       &
        \mapsto
      &
      \left(
        \begin{array}{cccc}
        0 & 0 & 0 & 1 \\
        1 & 1 & 1 & 0 \\
        1 & 1 & 1 & 0 \\
        1 & 1 & 1 & 0 \\        
        \end{array}
      \right)
      &
        \mapsto              
      \left(
        \begin{array}{cccc}
        1 & 1 & 1 & 1 \\
        1 & 1 & 1 & 1 \\
        1 & 1 & 1 & 1 \\
        1 & 1 & 1 & 1 
        \end{array}
      \right)
    \end{array}.   
\]
\par
\textbf{Right merging:}
\[
    \begin{array}{ccccccc}
      \left(
        \begin{array}{cccc}
        0 & 0 & 0 & 1 \\
        0 & 1 & 0 & 0 \\
        0 & 0 & 1 & 0 \\
        1 & 0 & 0 & 0 
        \end{array}
      \right)
      &
        \mapsto
      &
      \left(
        \begin{array}{cccc}
        0 & 0 & 0 & 1 \\
        0 & 1 & 1 & 0 \\
        0 & 1 & 1 & 0 \\
        1 & 0 & 0 & 0         
        \end{array}
      \right)
       &
        \mapsto
      &
      \left(
        \begin{array}{cccc}
        0 & 1 & 1 & 1 \\
        0 & 1 & 1 & 1 \\
        0 & 1 & 1 & 1 \\
        1 & 0 & 1 & 0 \\        
        \end{array}
      \right)
      &
        \mapsto
              
      \left(
        \begin{array}{cccc}
        1 & 1 & 1 & 1 \\
        1 & 1 & 1 & 1 \\
        1 & 1 & 1 & 1 \\
        1 & 1 & 1 & 1 
        \end{array}
      \right)
    \end{array}.
\]
\end{example}
 
\subsection{Bracketing of permutations}\label{subsec:Bracketing}

Given a permutation $\pi= a_1a_2\dots a_n$ of $[n]$, Shapiro and Stephens, in \cite{SS}, show how to determine if $\pi$ is \full, or more generally, to determine the final configuration of bootstrap percolation, by bracketing the string ``$a_1 a_2 \dots a_n$'' while avoiding the need to compute the sequence of matrices (\ref{eqn:stringOfMs}).  In effect, the bracketing keeps track of the progress of a \tile\ merging algorithm by encoding tiles as substrings of $a_1 a_2 \dots a_n$. The key idea is given in Lemma \ref{lem:SSBracketingIdea}, for which the following definition is useful:

\begin{definition}\label{def:T_tRem}
 Let $\pi$ be a permutation of $[n]$, let $M_0$ be the permutation matrix of $\pi$, and let $T_0$ be a square \tile\ of size $k$ ($1\leq k\leq n$) in $M_0$. Consider a complete sequence of percolation steps as in (\ref{eqn:stringOfMs}).  For each $0 < t \leq r$, let $T_t$ be the \tile\ in $M_t$ occupying the same positions as $T_0$ occupies in $M_0$. We call $T_t$ the $t^{\text{th}}$ \debi{descendant} of $T_0$ and $T_0$ the \debi{origin} of $T_t$.  Note that if $T_0$ contains $k$ non-zero entries of $M_0$, then $T_0$ can be viewed as a permutation matrix for a permutation of $[k]$, since each row and column of $M_0$ (and hence $T_0$) contains only one non-zero entry.  If in addition $T_0$ represents a \underline{\full} permutation of $k$, then for some $t$ we will have that $T_t$ is a unitary square \tile\ of size $k$; in this case we call $T_t$ a \debi{good} \tile.
\end{definition}

\begin{example}\label{exmp:ExamplesOfTiles}
The following example shows the complete percolation of a $7\times 7$ permutation matrix.  The final configuration (to the right) consists of five square unitary \tiles, four of size $1$ and one of size $3$.  All five of these tiles are good (in agreement with the following Lemma), but the $4\times 4$ \tile\ in the northwest, comprising the cells $c_{i,j}$ such that $1\leq i, j \leq 4$, although it represents a permutation of of $[4]$, is not a good tile, since the permutation it represents is not \full.

\begin{figure}[h!]
\centering
\includegraphics[width=100mm]{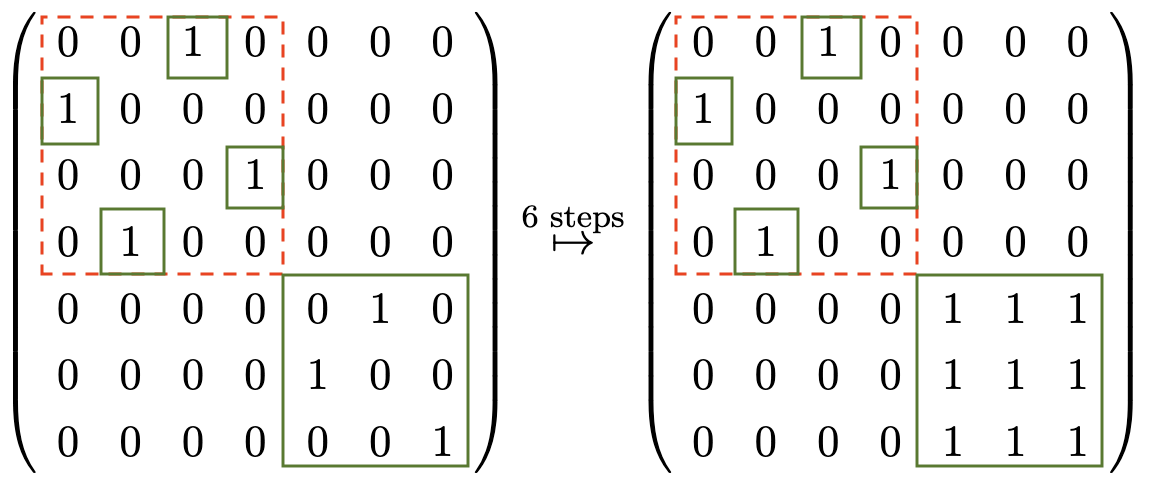}
\end{figure}

\end{example}

\begin{lemma}\label{lem:SSBracketingIdea}
    Let $\pi = a_1 a_2 \dots a_n$ be a permutation of $[n]$ and let $M_0$ be the permutation matrix of $\pi$.  If a \tile\ merging algorithm is run on $M_0$, then every square unitary \tile\ $T$ that appears during the computation is a good tile.  
    \end{lemma}
\begin{proof}
      The Lemma is clearly valid for the $n$ square unitary \tiles\ of size $1$ that are the unique nonzero cells of $M_0$.  Arguing by contradiction, we suppose the Lemma is false and consider a square unitary \tile\ $T$ of minimum size $k>1$ such that: 1) $T$ results from the merging of two previously-appearing square unitary \tiles\ $T' \text{ and } T''$, and 2) $T$ is not a good \tile.  By the minimality of $k$, we know that $T' \text{ and } T''$ have sizes $k' < k \text{ and } k'' < k$, respectively, hence are good \tiles.  Consequently, their origins $T'_0 \text{ and } T''_0$ contain $k'$ and $k''$ nonzero entries, respectively, so the origin $T_0$ of $T$ contains $k'+k'' = k$ nonzero entries and therefore represents a permutation of $[k]$.  Moreover, since $T$ is obtained from the merging of $T' \text{ and } T''$, it is clear that $T_0$ represents a \full\ permutation; consequently, $T$ is a good \tile, contradicting the hypothesis and completing the proof of the Lemma.    
\end{proof}

\begin{remark}\label{rem:meldFacts}
    Lemma \ref{lem:SSBracketingIdea} implies that every square unitary \tile\ $T$ of size $k$ that appears during the execution of a \tile\ merging algorithm corresponds uniquely to a length-$k$ substring ``$a_{j} a_{j+1} \dots a_{j+k-1}$'' that encodes the list of nonzero cells in $T_0$.  For instance, in Example \ref{exmp:leftRightBlockMerge}, the \tile\ of size $2$ that results from the first step of both left and right merging corresponds to the substring ``$32$'' of $\pi=1324$, whereas the substrings corresponding to the \tiles\ of size $3$ resulting from the second steps differ: for left merging, that substring is ``$132$,'' and for right merging, it is ``$324$.'' This is the key idea behind the bracketing process. We will call the (suitably bracketed) substrings representing the good \tiles\ ``\melds'' (they are called ``blocks'' in \cite{SS}).  A formal definition follows:  
\end{remark}

\begin{definition}\label{def:MeldDef}
    Let $\pi = a_1 a_2 \dots a_n$ be a permutation of $[n]$.  We inductively define a \dbi{\meld} of $\pi$ as follows: 
    \begin{itemize}
      \item[i.] Every component $a_j$ of $\pi$ is a \meld, and
      \item[ii.] If $m_1 \text{ and } m_2$ are consecutive \melds\ such that the union of the components of $m_1 \text{ and } m_2$ is a (permuted) set of consecutive integers in $[n]$\footnote{This condition is equivalent to the \tiles\ corresponding to $m_1 \text{ and } m_2$ being diagonally adjacent.}, then $(m_1, m_2)$ (resp.\ $[m_1, m_2]$) is the \meld\ resulting from their merger in case $\max(m_1) = \min(m_2)-1$ (resp.\ $\min(m_1) = \max(m_2)+1$).
    \end{itemize}
\end{definition}

To complete the bracketing of $\pi$, one mimics the \tile\ merging algorithm being used, but instead of merging square unitary \tiles\ in the percolating matrix, one merges the \melds\ that encode those \tiles.  To illustrate this, we mimic the computations in Example \ref{exmp:leftRightBlockMerge} to compute the left and right bracketings of $\pi=1324$.
\[
  \begin{array}{cccccccc}
\text{Left bracketing:} & \text{1 3 2 4} & \mapsto & \text{1 [3 2] 4} & \mapsto & \text{(1 [3 2]) 4} & \mapsto & \text{((1 [3 2]) 4)} \vspace{.05in}\\
\text{Right bracketing:} & \text{1 3 2 4} & \mapsto & \text{1 [3 2] 4} & \mapsto & \text{1 ([3 2] 4)} & \mapsto & \text{(1 ([3 2] 4))} 
  \end{array} 
\]
In this case, since the bracketings terminate with a single \meld\ encompassing all four digits in $\pi$, we confirm our earlier observation that $\pi$ is \full.

\subsection{ ``Top-level bracketings'' of a \full\ permutation}\label{subsec:TopLevelBracketings}

\begin{lemma}\label{lem:topLevelBracketing}
    Let $\pi$ be a \full\ permutation of $[n]$, for $n>1$.  Then, for any bracketing algorithm, $\pi$ has a ``top-level'' bracketing of the form $(m_1, m_2)$ or $[m_1, m_2]$ that comprises the integers in $[n]$.
\end{lemma}

\begin{proof}
      The hypothesis implies that the bracketing process on  $\pi$ must terminate with a single \meld\ comprising the integers in $[n]$. The initial number of \melds\ is $n\geq 2$, corresponding to the $n$ digits in $\pi$, and this number decreases by 1 every time two \melds\ are merged.  Since the process terminates with a single \meld\ that must arise from the merging of two \melds\ $m_1 \text{ and }m_2$, the result follows.     
\end{proof}

\begin{remark}\label{rem:topLevelBracketingShapes}
    Given our convention on how a permutation $\pi$ of $[n]$ is mapped to its permutation matrix (see Remark \ref{rem:permToPermMatrix}), one checks that the orientation of the \tiles\ $T_1 \text{ and } T_2$ corresponding to the \melds\ $m_1 \text{ and }m_2$ in the top-level bracketing have the following orientations prior to their merger:
\begin{equation}
    \begin{array}{ccc}
      (m_1, m_2) & \leftrightarrow & \left(
                            \begin{array}{cc}
                                 0 & T_2\\
                                 T_1 & 0
                            \end{array}\right)\vspace{.03in} \\
     \text{[}m_1, m_2\text{]} & \leftrightarrow & \left( 
                            \begin{array}{cc}
                                 T_1 & 0\\
                                 0 & T_2
                            \end{array}\right).
    \end{array}
\end{equation}
\end{remark}

\subsection{The binary tree associated to the bracketing of a \full\ permutation}\label{subseq:BinaryTree}

As explained in \cite[p.\ 277]{SS}, a bracketing of a \full\ permutation can be displayed as a binary tree in which the vertices are the \melds\, and the left (resp.\ right) child of a meld of the form $(m_1, m_2)$ or $[m_1, m_2]$ is $m_1$ (resp.\ $m_2$), and the singleton melds are the leaves of the tree.  For example, here is the tree associated to the \full\ (left-bracketed) permutation $\pi = ((1[32])4)$:

\setlength{\unitlength}{1in}
\begin{picture}(6,1.7)
    \put(2.5,1.4){((1 [3 2]) 4)}
    \put(2.1,1){(1 [3 2])}
    \put(3.3,1){4}
    \put(1.9,.6){1}
    \put(2.5,.6){[3 2]}
    \put(2.2, .2){3}
    \put(2.9,.2){2}
    \put(2.6,1.3){\line(-1,-1){.15}}
    \put(3.1,1.3){\line(1,-1){.15}}
    \put(2.2,0.9){\line(-1,-1){.15}}
    \put(2.4,0.9){\line(1,-1){.15}}
    \put(2.7,.5){\line(1,-1){.15}}
    \put(2.5,.5){\line(-1,-1){.15}}
\end{picture}
 
\begin{lemma}\label{lem:RtChildFact}
    Let $n>1$ and let $\pi$ be a permutation of $n$.  If left bracketing is performed on $\pi$, the corresponding binary tree has the following property: the right child $m_2$ of a non-leaf node $(m_1, m_2)$ (resp.\  $[m_1, m_2]$) is either a singleton or has the form $[m_{21}, m_{22}]$ (resp.\ $(m_{21}, m_{22})$); in brief, the type of the right child of a non-leaf node under left bracketing differs from the type of its parent.  Symmetrically, the left child of a non-leaf node under right bracketing differs from the type of its parent.
\end{lemma}

\begin{proof}
    Suppose to the contrary that under left bracketing, a non-leaf node has the form $(m_1, m_2)$, where $m_2 = (m_{21}, m_{22})$.  Recall from Remark \ref{rem:meldFacts} that each of the \melds\ $m_1$, $m_{21}$, $m_{22}$, and $(m_1, m_2)$ comprises a (permuted) consecutive run of integers in $[n]$, and the bracketing indicates that the integers in $m_1$ are all less than the integers in $m_{2}$ and that the integers in $m_{21}$ are all less than those in $m_{22}$.  A moment's reflection then implies that (under left merging) $m_1$ and $m_{21}$ will be merged to form $(m_1, m_{21})$ \textit{before} $m_{21} \text{ and } m_{22}$ can be merged, contradicting the hypothesis.  The proof in the case of a supposed non-leaf node of the form $[m_1, m_2]$ with $m_2 = [m_{21}, m_{22}]$ is similar, as is the proof of the last statement in the case of right bracketing. 
\end{proof}

\begin{remark}\label{rem:SS"error"}
The left-bracketing version of the assertion in Lemma \ref{lem:RtChildFact} is stated at the top of page 278 of \cite{SS}; subsequently, it is needed there to establish the count of \full\ permutations of $n\geq 2$.  However, this result need not be true depending on one's interpretation of the (somewhat vague) bracketing procedure given on page 277 of \cite{SS}:
\begin{equation}\label{eqn:SSAlgorithm}
  \begin{array}{l}
    \text{The permutation $\pi$ can be considered as a sequence of [\melds].}\\
    \text{We read $\pi$ left to right and use [the merging rule (ii.\ of Definition \ref{def:MeldDef})]}\\
    \text{to form new [\melds] whenever possible.}
  \end{array}  
\end{equation}
The following example makes the key issue clear.
\end{remark}

\begin{example}\label{exmp:SSAmbiguity}
Consider bracketing the permutation $\pi=4231$.  Beginning at the left, we find that the first adjacent pair of \melds\ that can be merged is $(2, 3)$.  Merging them yields $4(23)1$, leaving two adjacent pairs of mergeable \melds\ $(4, (23))$ and $((23), 1)$.  There are now two possible ways to complete the bracketing:
\[
    \begin{array}{cccccc}
     \text{i.} & 4(23)1 & \mapsto & [4(23)]1 & \mapsto & [[4(23)]1],\\
     \text{ii.} & 4(23)1 & \mapsto & 4[(23)1] & \mapsto &[4[(23)1]].
    \end{array}
\]
Both computations are compatible with the procedure sketched in (\ref{eqn:SSAlgorithm}), but only the first (the result of left merging) satisfies the left-bracketing version of Lemma \ref{lem:RtChildFact}.  The second computation (in which newly-merged \melds\ are, if possible, merged with their right-hand neighbors as a left-to-right traversal of $\pi$ continues) appears, based on the example given on page 277 of \cite{SS}, to be what was intended by the procedure (\ref{eqn:SSAlgorithm}).
\end{example}

In the sequel, we will have occasion to use the left-bracketing version of Lemma \ref{lem:RtChildFact}, and will only do so when left bracketing has been used.

\section{Indecomposable permutations}\label{sec:IndecPerms}

In this section we recall the definition and basic theory of indecomposable permutations.
The concept is due to Comtet in the early 1970's (See \cite{C1} or \cite{GKZ}).  Indecomposable permutations are sometimes 
known as \textit{irreducible} or \textit{connected} permutations. 

\subsection{Definition of indecomposable permutation}\label{subseq:IndecPermDef}

\begin{definition} A permutation $\pi \text{ of } [n]$ is called \debi{indecomposable} if there does not exist an index $1 \leq k < n$ such that $\pi(\{1, 2, \dots, k\}) = \{1, 2, \dots, k\}$. That is, there does not exist a proper initial segment of $\{1, 2, \dots, n\}$ that is mapped onto itself by $\sigma$; equivalently, no length-$k$ prefix of the single-line representation of $\pi$ is a permutation of $[k]$ for some $k<n$.  Otherwise we say $\pi$ is \debi{decomposable}. 
\end{definition}

\begin{example}\label{exmp:IndecAndDecPerms}
The permutation $\tau = 4132675$ is decomposable because $4132$ is a permutation of $[4]$, while the permutation $\sigma = 4167523$ is indecomposable.
\end{example}

Suppose $\pi = a_1 a_2 \dots a_n$ is a permutation of a set of integers (or any linearly ordered set). Following \cite{GKZ}, we define the
\textbf{\textit{reduced form}} of $\pi$ to be the permutation of $[n]$ obtained by replacing the $i^{\text{th}}$ smallest element of $\pi$ with $i$.  We denote the reduced form of $\pi$ by $\pi_{red}$. For example, if $\pi = 2754$, then $\pi_{red} = 1432$. This notion permits us to extend the definition of indecomposable permutation to permutations of arbitrary (finite) linearly ordered sets: given any such partition $\pi$, we say that $\pi$ is \textbf{\textit{indecomposable}} provided that $\pi_{red}$ is indecomposable.  
\begin{example}
   $\pi=2754$ is decomposable because its reduced form $1432$ is decomposable.
\end{example}

\subsection{Indecomposable components of a permutation}\label{subsec:DefOfIndecComps}

Let $\pi= a_1 a_2 \dots a_n$ be a permutation of a linearly ordered set, and let $\pi_{red} = a'_1 a'_2 \dots a'_n$. Then let $k$ denote the least positive integer such that $a'_1 a'_2 \dots a'_k$ is a permutation of $[k]$, and let $\pi_1 = a_1 a_2\dots a_k$.  One then checks easily that:  
     \begin{enumerate}
       \item[i.] $(\pi_1)_{red} = a'_1 a'_2 \dots a'_k$,
       \item[ii.] $a'_1 a'_2 \dots a'_k \text{ and (hence) } \pi_1$ are indecomposable, 
       \item[iii.] $\pi_1$ is the maximum-length indecomposable prefix of $\pi$, and 
       \item[iv.] $\pi$ is indecomposable $\Leftrightarrow\ \pi=\pi_1$.
     \end{enumerate}

\newcommand{\comps}{\operatorname{comps}}
It is known (e.g., see \cite{GKZ}) that for any permutation $\pi$ (taking the product to be concatenation of the single line notation of the factors), 
we have
\[
    \pi = \pi_1\pi_2\cdot...\cdot\pi_q,
\]
for some $q\geq 1$, such that $\pi_1$ is the indecomposable prefix of $\pi$ of maximum length and $\pi_j$ is indecomposable for $j>1$. The permutations $\pi_j$ are called the \textbf{\textit{(indecomposable) components}} of $\pi$.  We can construct them inductively as follows: Given a permutation $\pi$, we let $\pi_1$ be the maximum-length indecomposable prefix of $\pi$, as above, and write $\pi= \pi_1 \pi'$.  Then we define
\[
    \comps(\pi) = \left\{ \begin{array}{l}
                            (\pi_1), \text{ if } \pi \text{ is indecomposable, and}\\
                                (\pi_1) \comps(\pi'), \text{ otherwise.}
                          \end{array}
                  \right.
\]
\begin{example}
    Let $\pi = 24135867$.  Then   
\[
  \comps(24135867) = (2413) \comps(5867) = (2413)(5)\comps(867) = (2413)(5)(867). 
\]
The components of $\pi$
are thus $\pi_1=(2413)$ (with $\pi_1 = {\pi_1}_{red}$), $\pi_2=(5)$ (with indecomposable reduced form $(1)$), and $\pi_3=(867)$ (with indecomposable reduced form $(312)$).
\end{example}

In section \ref{subsec:IndecCompOfFullPermutations}, we show how the indecomposable components of a \full\ permutation can be computed by using left bracketing.

\section{Full Indecomposable Permutations and the ``Half-lemma''}\label{sec:IndecCompOfFullPermutations}

The main goal of this section is to prove that exactly half of the \full\ permutations of $[n], \text{for } n \geq 2,$ are indecomposable.  We also show that the indecomposable components of a \full\ permutation can be computed by left bracketing, and prove that a permutation is \full\ if and only if its indecomposable components are \full. 

\subsection{An involution on \full\ permutations}\label{subsec:InvOnFullPerms}

Given a permutation $\pi = a_1 a_2 \text{...} a_n$, we define the \debi{reversal} of $\pi$ to be the permutation $\pi^R= a_n
a_{n-1} \text{...} a_1$. It is clear that the reversal map $R:\pi \mapsto \pi^R$ is an involution on the set of all permutations. Furthermore, it is easy to see that $R$ is an involution on the set of \full\ permutations:

\begin{lemma}\label{lem:ReverLem1} If $\pi$ is a \full\ permutation of $[n]$, then so is $\pi^R$, and conversely.
\end{lemma}

\begin{proof} By Theorem \ref{thm:InvarProp}, the final configuration resulting from bootstrap percolation is independent of the order
in which cells are mutated. That being so, it suffices to 
note that running the left merging algorithm on $\pi^R$ is the ``mirror image'' of running the right merging algorithm on $\pi$, so if the latter's final configuration is a square unitary \tile\ of size $n$, so is the former's, and vice versa. 
\end{proof}

To clarify the preceding proof, we offer the following 

\begin{example}\label{exmp:L-RandR-Lmerging}
Let $\pi = 31254$ and $\pi^R = 45213$. We now perform the first three steps of  right merging on $\pi$ and left 
merging on $\pi^R$.  After each step, the results are mirror images of one another.  In this case, both $\pi$ and $\pi^R$ are \full, since the unitary \tiles\ of sizes $2 \text{ and } 3$ in the rightmost matrices will merge to yield unitary \tiles\ of size $5$.
\par
\bigskip

Right merging on $\pi$: 
{\scriptsize
\[
  {
  \begin{array}{ccccccc}
    \left(\begin{array}{ccccc}
        0 & 0 & 0 & 1 & 0 \\
        0 & 0 & 0 & 0 & 1 \\
        1 & 0 & 0 & 0 & 0 \\
        0 & 0 & 1 & 0 & 0 \\
        0 & 1 & 0 & 0 & 0
    \end{array}\right) & \mapsto & 
    \left(\begin{array}{ccccc}
        0 & 0 & 0 & 1 & 1 \\
        0 & 0 & 0 & 1 & 1 \\
        1 & 0 & 0 & 0 & 0 \\
        0 & 0 & 1 & 0 & 0 \\
        0 & 1 & 0 & 0 & 0
    \end{array} \right) & \mapsto &
    \left(\begin{array}{ccccc}
        0 & 0 & 0 & 1 & 1 \\
        0 & 0 & 0 & 1 & 1 \\
        1 & 0 & 0 & 0 & 0 \\
        0 & 1 & 1 & 0 & 0 \\
        0 & 1 & 1 & 0 & 0
     \end{array}\right) & \mapsto &  
    \left( \begin{array}{ccccc}
        0 & 0 & 0 & 1 & 1 \\
        0 & 0 & 0 & 1 & 1 \\
        1 & 1 & 1 & 0 & 0 \\
        1 & 1 & 1 & 0 & 0 \\
        1 & 1 & 1 & 0 & 0
    \end{array} \right)
  \end{array} }
\]
}

Left merging on $\pi^R$:
{\scriptsize
\[
  {
  \begin{array}{ccccccc}
    \left(\begin{array}{ccccc}
        0 & 1 & 0 & 0 & 0 \\
        1 & 0 & 0 & 0 & 0 \\
        0 & 0 & 0 & 0 & 1 \\
        0 & 0 & 1 & 0 & 0 \\
        0 & 0 & 0 & 1 & 0 
        \end{array} \right)  & \mapsto & 
    \left( \begin{array}{ccccc}
        1 & 1 & 0 & 0 & 0 \\
        1 & 1 & 0 & 0 & 0 \\
        0 & 0 & 0 & 0 & 1 \\
        0 & 0 & 1 & 0 & 0 \\
        0 & 0 & 0 & 1 & 0
    \end{array} \right)  & \mapsto &
    \left(\begin{array}{ccccc}
        1 & 1 & 0 & 0 & 0 \\
        1 & 1 & 0 & 0 & 0 \\
        0 & 0 & 0 & 0 & 1 \\
        0 & 0 & 1 & 1 & 0 \\
        0 & 0 & 1 & 1 & 0
    \end{array} \right)  & \mapsto &
    \left(\begin{array}{ccccc}
        1 & 1 & 0 & 0 & 0 \\
        1 & 1 & 0 & 0 & 0 \\
        0 & 0 & 1 & 1 & 1 \\
        0 & 0 & 1 & 1 & 1 \\
        0 & 0 & 1 & 1 & 1
    \end{array} \right)
  \end{array} }
\]
}
\end{example}
\bigskip

Our next goal is Theorem \ref{th:ReversalTh}, which states that the restriction of the involution $R$ to \full\ permutations 
sends indecomposable permutations to decomposable permutations and vice versa.  For the proof we need the following

\begin{lemma}\label{lem:ReversalLemma} 
If a permutation $\pi$ of $[n], \text{ for } n\geq 2$, is \full, then $\pi$ is indecomposable if and only if the top-level
bracketing (see Lemma {\rm \ref{lem:topLevelBracketing}}) of $\pi$ has the form $\left[m_1,m_2\right]$. Equivalently, $\pi$ is decomposable if and only if the top-level bracketing of $\pi$ has the form $\left(m_1,m_2\right)$.
\end{lemma}

\begin{proof}
It suffices to prove one of the equivalent statements. We prove the second.
\par
$\bm{(\Rightarrow)}$ Let $\pi$ be a \full\ decomposable permutation of $[n ], \text{ for } n \geq 2$, and let $m_1$ and $m_2$
be the \melds\ in the top-level bracketing of $\pi$. By definition of decomposable, we know that $\pi$ maps $[k]$ to $[k]$ for some $1\leq k<n$.
We must show that the top-level bracketing has the form $\left(m_1, m_2\right).$ Suppose to the contrary that the top-level bracketing of $\pi$ has the form $\left[m_1, m_2\right]$ where $m_1$ is a \meld\ of length $p < n$ comprising the $p$ largest values in $[n]$; that is, $\pi$ maps $[p]$ to $\{n-p+1, n-p+2, \dots, n \}$. We now have the following cases:\\
{\bf  Case 1:} $k\leq p$. It follows that 
\[
    \{1, 2, \dots, k\} \subseteq \{n-p+1, n-p+2, \dots, n\} \Rightarrow 1= n-p+1 \Rightarrow p = n, \text{ contradiction.}
\]
{\bf  Case 2: }$k > p$. It follows that
\[
    \{n-p+1, n-p+2, \dots, n\} \subseteq \{1, 2, \dots, k\} \Rightarrow k = n, \text{ contradiction.}
\]
Since both cases lead to contradictions, we have that the top-level bracketing has the form $\left(m_1, m_2\right)$, as desired.

$\bm{(\Leftarrow )}$ Now suppose that $\pi$ is a \full\ permutation of $[n]$ such that its top-level bracketing
has the form $\left(m_1, m_2\right)$, where $m_1$ is a \meld\ of length $k < n$ representing a permutation of $[k]$. It follows by definition that $\pi$ is decomposable. This completes
the proof of the Lemma.
\end{proof}

\begin{remark}\label{rem:topLevBracketingInvariance}
We note the following:
  \begin{enumerate} 
     \item[i.] A moment's reflection shows that the results of left bracketing (or other bracketing algorithm) on a permutation $\pi$ of $\{k+1, k+2, \dots, k+n\}$ will be identical in structure to the results of the same bracketing algorithm on the reduced form $\pi_{red}$.  For example, let $\pi = 68745$.  Then
     \[
       \begin{array}{rcl}
     \text{left bracketing of } \pi & = & [(6[87])(45)], \text{ and}\\
     \text{left bracketing of } \pi_{red} & = & [(3[54])(12)].
       \end{array}
     \]
     In particular,  Lemma \ref{lem:ReversalLemma} holds for a \full\ permutation of $\pi = \{k+1, k+2, \dots, k+n\}$.
      \item[ii.]The proof of Lemma \ref{lem:ReversalLemma} did not reference the \tile\ merging or corresponding bracketing algorithm used. It follows that for \full\ permutations $\pi$, the form (\textit{i.e.}, $(m_1 , m_2)$ or $[m_1 , m_2]$) of $\pi$'s top-level bracketing is independent of the bracketing algorithm, although the \melds\ $m_1 \text{ and } m_2$ are not, as shown, for instance, by Example \ref{exmp:SSAmbiguity}.
\end{enumerate}
\end{remark}

\subsection{The Half-lemma}\label{subseq:halfLemmaSec}

\begin{theorem}\label{th:ReversalTh}For $n\geq  2$, the involution $R$, restricted to the set of \full\ permutations of $[n]$, carries indecomposable permutations to decomposable permutations and vice versa.
\end{theorem}
\begin{proof} Let $\pi$ be a \full\ indecomposable permutation of $[n], \text{ for } n \geq 2$. We must show that $\pi^R$ is decomposable. By Lemma \ref{lem:ReversalLemma}, the top-level bracketing of $\pi$ has the form $\left[m_1, m_2\right]$ for any bracketing algorithm (although the melds $m_1 \text{ and } m_2$ may differ with different algorithms).  Whence, by Remark \ref{rem:topLevelBracketingShapes}, the square unitary \tiles\ $T_1 \text{ and } T_2$ corresponding to $m_1 \text{ and } m_2$ that appear during right  merging have the form    
$\left(
\begin{array}{cc}
 T_1 & 0 \\
 0 & T_2 
\end{array}
\right)$ (prior to being merged). Then, as observed in the proof of Lemma \ref{lem:ReverLem1}, left merging on $\pi^R$
is the mirror image of right merging on $\pi$. Thus the top-level bracketing of
$\pi^R$ will correspond to be the mirror image of $\left(
\begin{array}{cc}
 T_1 & 0 \\
 0 & T_2 
\end{array}
\right)$, which is $\left(
\begin{array}{cc}
 0 & T_1 \\
 T_2 & 0 
\end{array}
\right)$. Accordingly, by Lemma \ref{lem:ReversalLemma}, $\pi ^R$ is decomposable. Similarly, if $\pi$ is decomposable, then $\pi ^R$ is indecomposable. {
}This completes the proof.
\end{proof}

As an immediate consequence, we obtain

\begin{corollary}[the Half-lemma]\label{cor:HalfLemma}  For $n\geq 2$, half of the \full\ permutations of $[n]$ are indecomposable and half are decomposable. \qed
\end{corollary}

\subsection{Indecomposable components of a \full\ permutation}\label{subsec:IndecCompOfFullPermutations}

Our goal in this section is to show that the indecomposable components of a \full\ permutation $\pi \text{ of } [n]$ can be obtained from the left bracketing of $\pi$.  To this end, we prove the following Lemma:

\begin{lemma}\label{lem:lastIndecComp}
    Let $\pi = a_1 a_2 \dots a_n$ be a permutation of $[n]$, and let $\pi = \pi_1 \pi_2 \dots \pi_q$ be the factorization of $\pi$ into indecomposable components. Then 
    \[
    \pi_q\ =\ \text{the maximum-length indecomposable suffix of } \pi.
    \]  
\end{lemma}

\begin{proof} 
  From the discussion in Section \ref{subsec:DefOfIndecComps}, one sees that $\pi_1$ is an indecomposable permutation of $[k]$ for some $k\leq n$; if $k<n$, one has that $\pi_2$ is an indecomposable permutation of $\{ k+1, \dots, k+l\}$ for some $l \leq n-k$, and so on.  In this way, we obtain that the rightmost indecomposable component $\pi_q = a_v a_{v+1} \dots a_n$ satisfies the following condition (where the LHS is taken to be $-\infty$ if $v=1$): 
    \[
       \max(\{a_1, a_2, \dots, a_{v-1}\}) < \min(\{a_v, a_{v+1},  \dots, a_n\}).
    \]
It is then clear that any suffix of $\pi$ strictly containing $\{a_v, a_{v+1},  \dots, a_n\}$ is necessarily decomposable, so we obtain that $\pi_q$ is the maximum-length indecomposable suffix of $\pi$.
\end{proof}

Lemma \ref{lem:lastIndecComp} leads to a way to compute the indecomposable components of a \full\ permutation $\pi$ of $[n>1]$ through left bracketing.

\begin{theorem}\label{thm:CompsThruBracketing}
    Let $\pi = a_1 a_2 \dots a_n$ be a \full\ permutation of $[n]$. If the top-level left bracketing of $\pi$ has the form $[m_1, m_2]$, then $\pi = \pi_1$ is indecomposable. Otherwise, $\pi$ is decomposable, the top-level left bracketing has the form $(m_1, m_2)$, and $m_2$ is the rightmost indecomposable component of $\pi$ (with its left bracketing).
\end{theorem}

\begin{proof}
  By Lemma \ref{lem:ReversalLemma} , we know that the top-level bracketing of a \full\ permutation of $[n>1]$ has the form $[m_1, m_2]$ (resp.\ $(m_1, m_2)$) if and only if $\pi$ is indecomposable (resp. decomposable); in particular, the first statement holds.  On the other hand, if $\pi$ is decomposable, the left bracketing has the form $(m_1, m_2)$, where $m_2$ is either a singleton or has the form $[m_{21}, m_{22}]$, by Lemma \ref{lem:RtChildFact}. In either case, we have that every element $a_j$ of $m_1$ is less than every element $a_l$ of $m_2$. Furthermore, $m_2$ is indecomposable: this is clear if $m_2$ is a singleton, and otherwise is a consequence of Lemma \ref{lem:ReversalLemma} and point {\rm i.} of Remark \ref{rem:topLevBracketingInvariance}; indeed, $m_2 = [m_{21}, m_{22}]$ is the bracketing of a \full\ permutation of $\{n-s+1, n-s+2, \dots, n\}$ for some $s<n$. Accordingly, $m_2$ is the bracketing of the maximum-length indecomposable suffix of $\pi$, so Lemma \ref{lem:lastIndecComp} yields that $m_2$ is the (bracketed) rightmost indecomposable component of $\pi$, as desired.
\end{proof}
    
\begin{remark}\label{rem:HowToFindIndecomFactors}
    Given a \full\ decomposable permutation $\pi$ of $n>1$, Theorem \ref{thm:CompsThruBracketing} enables us to extract $\pi$'s indecomposable components $\pi_1 \pi_2\dots \pi_q$ directly from its left bracketing, as follows: From the top-level bracketing $(m_1, m_2)$, we obtain $m_2 = \pi_q$, so $\pi= m_1 \pi_q$. If $m_1$ is indecomposable, then $\pi_{q-1} = \pi_1=m_1$; otherwise $m_1=(m_{11}, m_{12})$ and so $\pi_{q-1} = m_{12}$, and so on.
\end{remark}

\begin{example}\label{exmp:indecompComps}
    We illustrate Remark \ref{rem:HowToFindIndecomFactors} with the following example: $\pi = 312645798$.  We first apply left bracketing:
\[
  \begin{array}{ccccc}
     3\, 1\, 2\, 6\, 4\, 5\, 7\, 9\, 8 & \mapsto & 3\, (1\, 2)\, 6\, 4\, 5\, 7\, 9\, 8 & \mapsto & [3\, (1\, 2)]\, 6\, 4\, 5\, 7\, 9\, 8\\
{} & \mapsto & [3\, (1\, 2)]\, 6\, (4\, 5)\, 7\, 9\, 8 & \mapsto & [3\, (1\, 2)]\, [6\, (4\, 5)]\, 7\, 9\, 8\\
{} & \mapsto & ([3\, (1\, 2)]\, [6\, (4\, 5)])\, 7\, 9\, 8 & \mapsto & (([3\, (1\, 2)]\, [6\, (4\, 5)])\, 7)\, 9\, 8\\
 {} & \mapsto & (
 ([3\, (1\, 2)]\, [6\, (4\, 5)])\, 7)\, [9\, 8] & \mapsto & ((([3\, (1\, 2)]\, [6\, (4\, 5)])\, 7)\, [9\, 8])
  \end{array}.  
\]
Reading the computation backwards, we obtain:
\[
   \begin{array}{rcl}
       \text{Rightmost component of } ((([3\, (1\, 2)]\, [6\, (4\, 5)])\, 7)\, [9\, 8]) & = & [9\, 8],\\
       \text{Rightmost component of } (([3\, (1\, 2)]\, [6\, (4\, 5)])\, 7) & = & 7,\\
       \text{Rightmost component of } ([3\, (1\, 2)]\, [6\, (4\, 5)]) & = & [6\, (4\, 5)],\\
       \text{Rightmost component of } [3\, (1\, 2)] & = & [3\, (1\, 2)] \text{ (it's indecomposable).}
   \end{array}
\]
It follows that the indecomposable factorization of $\pi$ is $(312)(645)(7)(98)$.  
\end{example}

The next result completes the agenda for this section of the paper.

\begin{theorem}\label{thm:piFullIFFCompsFull}
    Let $\pi$ be a permutation of $[n]$, for $n\geq 1$.  Then $\pi$ is \full\ if and only if each of its indecomposable components is \full.
\end{theorem}

\begin{proof}
  Let $\pi$ be a permutation and $\pi_1 \pi_2 \dots \pi_q$, for $q\geq 1$, its factorization into indecomposable components.  If $q=1$ (that is, $\pi$ is indecomposable), the theorem is immediate, so it suffices to show the result holds if $q>1$ ($\Leftrightarrow$ $\pi$ is decomposable).  Let   
  $m_j$ denote the left-bracketing of $\pi_j$, for $1\leq j\leq q$.  
  
  $\bm{(\Rightarrow)}$ Suppose $\pi$ is a \full\ decomposable permutation of $[n]$.  Then by Theorem \ref{thm:CompsThruBracketing} and Remark \ref{rem:HowToFindIndecomFactors}, we have that the left bracketing of $\pi$ has the form 
\[
    ((\dots((m_1, m_2), m_3),\dots, m_{q-1}),  m_q). 
\]
Recall that the left bracketing is a record of the progress of left merging beginning with the permutation matrix $M_0$ of $\pi$.
In particular, each of the \melds\ $m_j$ represents a \full\ permutation corresponding to a square unitary \tile, that is, each $\pi_k$ is \full, as desired. 

$\bm{(\Leftarrow)}$ Suppose $\pi$ is a decomposable permutation of $[n]$ with factorization into indecomposable components $\pi_1 \pi_2\dots \pi_q$, for $q>1$, and such that each factor $\pi_j$ is a \full\ permutation.  Let $m_j$ denote the left bracketing of $\pi_j$ for $1\leq j \leq q$. From the discussion at the start of the proof of Lemma \ref{lem:lastIndecComp}, we know that there is an increasing sequence of integers 
\[
    1\leq s_1 < s_2 < \dots< s_{q-1} < n
\]
such that
\[
  \begin{array}{c}
    m_1 \text{ is a bracketed permutation of } \{1, 2, \dots, s_1\},\\
    m_2 \text{ is a bracketed permutation of } \{s_1+1, s_1+2, \dots, s_2\},\\
    \vdots\\
    m_{q-1} \text{ is a bracketed permutation of } \{s_{q-2}+1, s_{q-2+2}, \dots, s_{q-1}\},\\
    m_q \text{ is a bracketed permutation of } \{s_{q-1}+1, s_{q-1}+2, \dots, n\}.
  \end{array}
\]
One can now carry out bootstrap percolation on the permutation matrix asociated to $\pi$ by first percolating the \tiles\ associated to the $\pi_j$; since the latter are by hypothesis \full\ permutations, one obtains a matrix containing a sequence of square unitary \tiles\ that are diagonally adjacent along the antidiagonal, which when merged yield a single square unitary \tile\ of size $n$, demonstrating that $\pi$ is \full.  This completes the proof.
\end{proof}

\begin{example}\label{exmp:exmp:indecompCompsCont}
    To illustrate the second part of the preceding proof, consider the permutation $\pi= (312)(645)(7)(98)$ of Example \ref{exmp:indecompComps}. One sees by inspection that the (parenthesized) indecomposable components are \full.  The permutation matrix of $\pi$ is shown to the left and the matrix following the percolation of the indecomposable components to the right.

\begin{center}
\includegraphics[width=112mm]{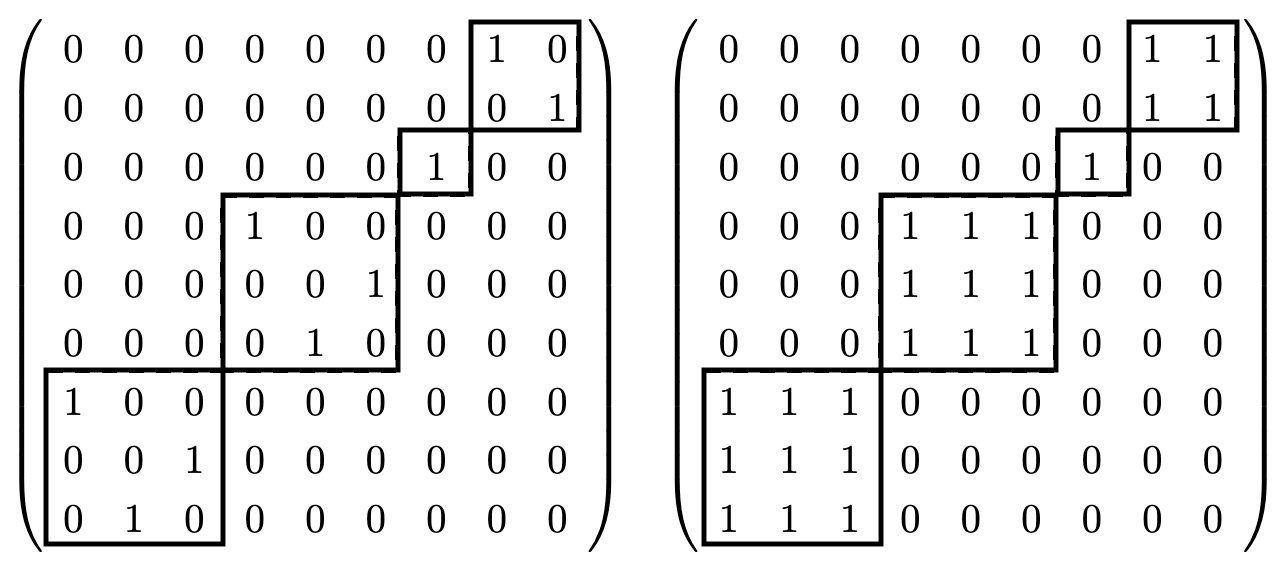}
\end{center}
\end{example}

\section{Counting \full\ permutations}\label{sec:CountFullPerms}

Let $I_n$ be the number of indecomposable permutations of $[n]=\{1,2,...,n\}$, let $p_n$ be the number of \full\ permutations of $[n]$, and let $q_n$ be the number of \full\ indecomposable permutations of $[n]$.  The goal of this section is to give a new proof of the following theorem, the main result of \cite{SS}:
\begin{theorem}\label{thm:MainThOfSS}
    The \full\ permutation matrices are enumerated by the sequence of (shifted) large Schr\"{o}der numbers.  That is, $p_n = S_{n-1}$.
\end{theorem}
The key ingredient for our proof is the following consequence of Corollary \ref{cor:HalfLemma}:
\begin{equation}\label{eqn:MainLemma}
    \begin{array}{l}
       \text{\textit{We have $p_1 = q_1 = 1$, and for $n>1$, we have $q_n = \frac{p_n}{2}$.}}
    \end{array}
\end{equation}

\begin{proof}
We begin by deriving a recursion for $p_{n+2}$. Let $\sigma$ be a \full\ permutation of $[n+2]$.  If $\sigma$ is indecomposable, we know there are $q_{n+2}$ choices for $\sigma$.  Otherwise, if $\sigma$ is decomposable, write $\sigma = \mu\tau$, where $\mu$ is the first indecomposable factor and $\tau$ is the rest of the permutation. 
We know $\sigma$ is \full\ if and only if both $\mu$ and $\tau$ are \full; this is an easy consequence of Theorem \ref{thm:piFullIFFCompsFull}. If $\mu$ is a permutation of $[k]$, there are $q_k$ 
choices for $\mu$ and $p_{n+2-k}$ choices for $\tau$.  Furthermore, when $\sigma$ is decomposable, $k$ must range between $1$ and 
$n+1$, inclusive. Summing all the choices, we have the following recursion:

\begin{equation} \label{VS1}
        p_{n+2} = q_{n+2} + \sum_{k=1}^{n+1} q_k \cdot p_{n+2-k} 
\end{equation}    
Now applying (\ref{eqn:MainLemma}), and breaking the term $k=1$ off from the sum, this becomes:
\begin{equation} \label{VS2}
        p_{n+2} = \frac{p_{n+2}}{2} + p_1 \cdot p_{n+1}  + \sum_{k=2}^{n+1} \frac{p_k}{2} \cdot p_{n+2-k} 
\end{equation} 
Multiply through by $2$ and use $p_1=1$ to obtain:
\begin{equation} \label{VS3}
        2p_{n+2} = p_{n+2} + 2 p_{n+1}  + \sum_{k=2}^{n+1} p_k \cdot p_{n+2-k} 
\end{equation} 
Thus,
\begin{equation} \label{VS4}
        p_{n+2} = 2 p_{n+1}  + \sum_{k=2}^{n+1} p_k \cdot p_{n+2-k} 
\end{equation} 
We can absorb one of the two terms equal to $p_{n+1}$ back into the sum as the $k=1$ term:
\begin{equation} \label{VS5}
        p_{n+2} = p_{n+1}  + \sum_{k=1}^{n+1} p_k \cdot p_{n+2-k} 
\end{equation} 
Now re-index the sum with $j = k-1$, which ranges from $0$ to $n$:
\begin{equation} \label{VS6}
        p_{n+2} = p_{n+1}  + \sum_{j=0}^{n} p_{j+1} \cdot p_{n+1-j} 
\end{equation}

Let's see what this recursion says for the shifted sequence $c_{n-1} = p_{n}$.  Rewriting in terms of $c_n$ our recursion now reads:

\begin{equation} \label{VS}
        c_{n+1} = c_{n}  + \sum_{j=0}^{n} c_{j} \cdot c_{n-j} 
\end{equation}

This is exactly the recursion defining the (large) Schr\"{o}der numbers $\{S_n\}$ \cite[p. 259]{West}.  Furthermore, $1 = S_0 = c_0 = p_1$.  This completes the proof of the theorem.
\end{proof}

The little Schr\"{o}der numbers (also called super-Catalan numbers or small Schr\"{o}der numbers): 1, 1, 3, 11, 45, 197, 903, 4279, 20793, \dots\ count, among other objects, the number of subdiagonal lattice paths from $(0, 0)$ to $(n, n)$ consisting of steps East $(1, 0)$, North $(0, 1)$ and Northeast $(1, 1)$ with no Northeast steps along the main diagonal. See \cite[A001003]{oeis} for more information about this sequence.

Since the little Schr\"{o}der numbers are exactly half the value of the large Schr\"{o}der numbers (the initial term excepted), applying statement (\ref{eqn:MainLemma}) again, we obtain:

\begin{corollary}
    The 
    \full\ indecomposable matrices are enumerated by the sequence of (shifted) little Schr\"{o}der numbers.
    That is, $q_n = s_{n-1}$.
\end{corollary}

\section{The \texorpdfstring{$n$}{n}-Kings Problem}\label{sec:NKingsProb}
The $n$-kings problem asks for the number of ways one can place $n$ kings on an $n \times n$ chessboard, with one king in
each row and column, so that none of the kings may attack each other.  If you think of the chessboard as an $n \times n$
matrix, with a $1$ entry where the kings are placed and other entries $0$, then because of the condition of one king
in each row and column, one obtains a permutation matrix. Clearly, none of the kings can attack any others if none of them are
diagonally adjacent (we already know they are not horizontally or vertically adjacent since they are in distinct rows and
columns), which means there are no mutable cells if bootstrap percolation is applied.  Thus, the permutation matrix is in 
a no-growth configuration if and only if the configuration is a solution to the $n$-kings problem.  Thus, as we stated in 
the Introduction, if $a_n$ is the number of solutions to the $n$-kings problem, then $a_n$ is also the number of 
$n \times n$ permutation matrices that are in a `no-growth configuration' with respect to bootstrap percolation. The 
goal of this section is to determine $a_n$ explicitly. 
\par
\bigskip
First, following \cite{SS}, we derive the generating function $A(t)$ of the sequence $\{a_n\}$ as the composition (\ref{eqn:keyCompForAn}) of two known functions. 
Our argument depends on the main result in \cite{V}, which shows how to express the Taylor coefficients $T_n(f\circ g;a)$ of 
a composite function $f\circ g$ (about $t=a$) in terms of the Taylor coefficients of the `factors.'  (All our series are 
powers of $t$, so $a=0$.) Second, we will apply this same result of \cite{V} to obtain an explicit formula for $a_n$ from 
this generating function, which is a new result.
\par
\bigskip
One of the known functions involved in our argument is the generating function $\eps(t)=\sum_{n \geq 0}n!t^n$ of the 
sequence $\{n!\}$, a power series with radius of convergence $0$. So, arguments that treat $\eps(t)$ as a function 
might be invalid --- this generating function only makes sense as a formal power series.  The arguments in \cite{V} 
relied on the function-theoretic viewpoint. However, we now have a proof of the main result of \cite{V} based 
entirely on formal power series (see \cite{V2}), so the calculations we are about to do are valid.\par 
\bigskip
To prove (\ref{eqn:keyCompForAn}), we'll count the $n \times n$  permutation matrices in two ways. On the one hand, we know there are $n!$ such matrices. On the other hand, we can count them by their final configurations after percolation. That is, we'll compute the number of permutations that percolate to each possible final configuration, and add these counts over all such configurations.
\par
\bigskip
We first consider how to count the number of permutations that percolate to a given final configuration.   
Recall from \S \ref{sec:finalConfig} that the final configuration $M_r$ of any permutation matrix following percolation consists of $m$ square unitary \tiles\ $\{T_1, T_2, \dots , T_m\}$, for $1\leq m \leq n$, in a no-growth configuration, and where the size of $T_j$ is $s_j$. Fix for the moment the no-growth configuration occupied by the $\{T_j\}$.  Then the possible final positions having $m$ square unitary \tiles\ occupying the given no-growth configuration are in one-to-one correspondence with the elements $\pi =(s_1,s_2,...,s_m) \in \mathcal{C}_{n,m}$ (the set of compositions of $n$ with $m$ parts) that specify the sizes of the \tiles\ from left to right. 
\par
\bigskip
We now determine the number of $n\times n$ permutation matrices that will percolate to the final configuration corresponding to the composition $\pi$. 
Since the percolation takes place independently within each unitary \tile\ (see \cite{V2}), it follows that the number of $n \times n$ permutation matrices which percolate to this particular $M_r$ is the 
product of the number that percolate to each unitary \tile\ $T_j$. Indeed, we know by Lemma \ref{lem:SSBracketingIdea} that each 
$T_j$ is a ``good'' \tile, meaning that it results from the percolation of a \full\ permutation of $[s_j]$, the number of which we denote by $p_{s_j}$ as in the previous section.  It follows that the number of initial matrices that percolate to the final configuration corresponding to $\pi$ is the product $\prod_{j=1}^m p_{s_j}$. Denoting the number of parts $s_j$ of $\pi$ that are equal to $i$ by $\pi_i$, we can rewrite the previous product as $\prod_{i=1}^n(p_i)^{\pi_i}$. 
\par
\bigskip
Summing the previous expression over all the $\pi\in \mathcal{C}_{n,m}$, we find that the number of  permutations that percolate to any final result having the current no-growth configuration is given by $\sum_{\pi \in \mathcal{C}_{n,m}} \prod_{i=1}^n (p_i)^{\pi_i}$.  Hence, multiplying by the number of no-growth configurations with $m$ \tiles, we obtain that the total number of permutations that percolate to any $m$-\tile{d} result is given by 
\[
a_m\cdot\sum_{\pi \in \mathcal{C}_{n,m}} \prod_{i=1}^n (p_i)^{\pi_i}.
\]
Finally, summing over $m$, we find that the number of permutations of $[n]$ is given by
\begin{equation} \label{countingpermutationmatrices}
    n! = \sum_{m=1}^n a_m \sum_{\pi \in \mathcal{C}_{n,m}} \prod_{i=1}^n (p_i)^{\pi_i}.
\end{equation}

Recall $\eps(t) = \sum_{n \geq 0} n!t^n$ is the ordinary generating function of the sequence $\{n!\}$, and $A(t)$ is 
the generating function for $\{a_n\}$.  If we let $B(t)$ stand for the ordinary generating function of $\{p_n\}$, 
then applying the main result of \cite{V}  (which expresses the $n^{\text{th}}$ Taylor coefficient $T_n(A\circ B;0)$ of 
the composition $A\circ B$ on the LHS in terms of the Taylor coefficients $T_m(A;0) = a_m$ of $A$ and 
$T_i(B;0) = p_i$ of $B$ on the RHS), 
we obtain:

\[ T_n(A\circ B;0) = \sum_{m=1}^n a_m \sum_{\pi \in \mathcal{C}_{n,m}} \prod_{i=1}^n (p_i)^{\pi_i}.
\]
 Equation (\ref{countingpermutationmatrices}) then yields that the $n^{\text{th}}$ Taylor coefficient of $A\circ B$ is $n!$,  proving that $A(B(t)) = \eps(t)$.

This argument is essentially an expansion of the heuristic argument outlined in Section 3 of \cite{SS} to obtain the
composition, but we think the result in \cite{V} makes the argument clear and precise. 
\par
\bigskip
Next, what is the sequence $\{p_n\}$?  We saw in the last section that the number of full $n \times n$ permutation matrices is 
$p_n = S_{n-1}$, the $(n-1)$'st (large) Schr\"oder number. Thus, $B(t) = tR(t)$, where $R(t)$ is the ordinary generating 
function of the Schr\"oder numbers (denoted $R(x)$ in \cite{SS}).  Combined with the above, we have:
\par
\bigskip
\begin{equation}\label{AB=E} 
    A(tR(t)) = \eps(t) 
\end{equation}
\par
\bigskip
We know that $B(t)=tR(t)=\frac{1-t-\sqrt{1-6t+t^2}}{2}$. Define 
$g(t)=\frac{t(1-t)}{1+t}$. We leave it to the reader to verify that $g(t)$ and $B(t)$ are inverse functions.  
It follows by composing both sides of (\ref{AB=E}) by $g$ on the right that:
\[A\circ B = \eps\]
\[A \circ B \circ g = \eps \circ g\]
Since $B \circ g$ is the identity function, we conclude $A = \eps  \circ g$:
\begin{equation} \label{A=Eg}
    A(t) = \eps(g(t))= \eps\biggl( \frac{t(1-t)}{1+t} \biggr)   
\end{equation}
This is the derivation of $A(t)$ in \cite{SS}, giving the generating function for the solution to 
the $n$-kings problem.
\par
\bigskip
We now aim to extract our explicit solution from this. Observe that $g(0)=0$, and by definition
we have that
\[T_m(\eps;0)=m!\]
As for the inner function, by the linearity of the operator $T_i(-;0)$, we obtain, for $i \geq 2$:
\[ T_i(g;0) = T_i\biggl( \frac{t(1-t)}{1+t};0 \biggr)\]
\[= T_i\biggl( \frac{t}{1+t};0 \biggr) - T_i\biggl( \frac{t^2}{1+t};0 \biggr) \]
Now, multiplying a function by powers of $t$ just shifts the coefficients of its series, so we can rewrite this as a 
linear combination of coefficients of an easily recognized geometric series:
\[ = T_{i-1}\biggl( \frac{1}{1+t};0 \biggr) - T_{i-2}\biggl( \frac{1}{1+t};0 \biggr)\]
\[ = (-1)^{i-1} - (-1)^{i-2}\]
\[= (-1)^{i-2}(-1 -1) = (-1)^{i-2}(-2) = 2(-1)^{i-1}\]
\par
\bigskip
For $i=1$, we have $T_1(g;0)=T_1\biggl( \frac{t}{1+t};0 \biggr)=1$, since $T_1\biggl( \frac{t^2}{1+t};0 \biggr)=0$.  
And for $i=0$ we have already observed that $g(0)=0$. Summarizing, we have:
\par
\bigskip
\[T_i(g;0) =\Biggl\{\substack{0  \hspace{.3 cm}\text{ if }i = 0\\ \\ 1 \hspace{.3 cm} \text{ if } i = 1\\ \\ 2(-1)^{i-1} \text{ if } i \geq 2} \]
\par
\bigskip
Applying the main result of \cite{V}, we obtain:
\par
\bigskip
\[a_n = T_n(\eps\circ g;0) = \sum_{m=1}^n T_m(\eps;0) \sum_{\pi \in \mathcal{C}_{n,m}} \prod_{i=1}^n(T_i(g;0))^{\pi_i}\]
\[ = \sum_{m=1}^n m! \sum_{\pi \in \mathcal{C}_{n,m}} \prod_{i=2}^n(2(-1)^{i-1}))^{\pi_i}\]
\par
\bigskip
Note that for any composition $\pi$, we have $T_1(g;0)^{\pi_1} = (1)^{\pi_1} = 1$, explaining why we can ignore the term
$i=1$ in the product.  Also, observe that
\[\prod_{i=2}^n(2(-1)^{i-1}))^{\pi_i} = \prod_{i=2}^n2^{\pi_i} \prod_{i=2}^n(-1)^{(i-1){\pi_i}} \]
\[ = 2^{\sum_{i \geq 2} \pi_i} \cdot (-1)^{\sum_{i \geq 2}i\pi_i - \sum_{i \geq 2} \pi_i}\]
\[ = 2^{\ell(\pi)-\pi_1}\cdot (-1)^{n-\pi_1} \cdot (-1)^{-\ell(\pi)+\pi_1}\]
\[ = 2^{\ell(\pi) - \pi_1} \cdot (-1)^{n-\ell(\pi)} = (-1)^n \frac{(-2)^{\ell(\pi)}}{2^{\pi_1}}\]
\par
\bigskip
It follows that
\par
\bigskip
\[a_n = \sum_{m=1}^n m! \sum_{\pi \in \mathcal{C}_{n,m}} (-1)^n \frac{(-2)^{m}}{2^{\pi_1}}\]
\par
\bigskip
Pulling the constants out of the sums, we obtain our explicit formula (first observed in 2016):
\par
\bigskip
\begin{theorem} \label{Explicitnkings}
    Let $a_n$ be the number of solutions to the $n$-kings problem (that is, the number of `no growth' $n \times n$ 
    permutation matrices). Then we have the explicit formula:
    \[a_n \hspace{.2 cm}=\hspace{.1 cm} (-1)^n \sum_{m=1}^n m!(-2)^m \sum_{\pi \in \mathcal{C}_{n,m}} \frac{1}{2^{\pi_1}}\] 
    \[ \hspace{.3 cm}=\hspace{.1 cm} (-1)^n \sum_{\pi \in \mathcal{C}_n} \frac{\ell(\pi)!(-2)^{\ell(\pi)}}{2^{\pi_1}} \]
\end{theorem}
\par
\bigskip
We illustrate the theorem for $n=5$, using the first expression. We first list the $16$ compositions of $5$, grouped by length:
\[ (5) \]
\[ (1,4),(2,3),(3,2),(4,1)\]
\[ (1,1,3),(1,3,1),(3,1,1),(1,2,2),(2,1,2),(2,2,1)\]
\[ (1,1,1,2),(1,1,2,1), (1,2,1,1),(2,1,1,1)\]
\[ (1,1,1,1,1)\]
Using this, we can easily compute the inner sum.  When $m=1$, we're on the first row, with only one composition $(5)$ for which
$\pi_1=0$, so the inner sum in this case is $1$. When $m=2$ (second row), we obtain $2\frac{1}{2^1} + 2\frac{1}{2^0} = 3$. In 
the same way for $m=3$ (third row), we obtain $\frac{9}{4}$ for the inner sum, and for $m=4$ we obtain $\frac{1}{2}$ and 
finally when $m=5$ we obtain $\frac{1}{32}$.
\par
\bigskip
Therefore, the theorem gives:

\[ a_5 = (-1)^5\biggl[1!(-2)^1\cdot 1 + 2!(-2)^2 \cdot 3 + 3!(-2)^3\cdot \frac{9}{4} + 4!(-2)^4\cdot \frac{1}{2} + 5!(-2)^5 \cdot \frac{1}{32}\biggr] \]
\[ = 2 - 24 + 108 - 192 + 120 = 14,\]
as expected.

\begin{remark}
    Equation (\ref{A=Eg}) yields a second explicit formula for $a_n$.  The idea is to expand 
\[
  \begin{array}{rcl}
     A(t) = \eps(g(t)) & = & \eps\biggl( \frac{t(1-t)}{1+t} \biggr)\\
            {}         & = & \sum_{m=0}^{\infty}m!\, t^m\,(1-t)^m\,(1+t)^{-m}\\
            {}         & = & \sum_{m=0}^{\infty}\biggl(m!\, t^m\, \left(\sum_{q=0}^{m} (-1)^q \binom{m}{q} t^q \right)\, 
            \left( \sum_{k=0}^{\infty} \binom{-m}{k} t^k \right) \biggr)
  \end{array}
\]
using the binomial theorem (for both non-negative and negative exponents), and then to extract the coefficient of $t^n$, which is $a_n$. (We thank Ira Gessel for suggesting this approach.) One obtains
\[
  \begin{array}{rcl}
        a_n & = & \sum_{m=0}^n m!\, \left( \sum_{q=0}^{n-m}\, (-1)^q\, \binom{m}{q}\, \binom{-m}{n-m-q} \right)\\
        {} & = & \sum_{m=0}^n m!\, \left( \sum_{q=0}^{n-m}\, (-1)^q\, \binom{m}{q}\, (-1)^{n-m-q}\, \binom{n-m-q+m-1}{n-m-q} \right)\\
        {} & = & \sum_{m=0}^n m!\, \left( \sum_{q=0}^{n-m}\, (-1)^{n-m}\, \binom{m}{q}\, \binom{n-q-1}{n-m-q} \right).
  \end{array}\footnote{Identity: $\binom{-r}{n} = (-1)^n \binom{n+r-1}{n}$}
\]
Rewriting the last expression by summing over $d=n-m$ instead of $m$, we obtain
\[
    a_n = \sum_{d=0}^n (n-d)!\, \left( \sum_{q=0}^{d}\, (-1)^{d}\, \binom{n-d}{q}\, \binom{n-q-1}{d-q} \right);
\]
replacing $d$ by $k$ and $q$ by $i$ in the last expression, one obtains the formula for $a_n$ found by Abramson and Moser \cite{AbramsonMoser1966} using a different approach; note that their outer sum runs from $0 \text{ to } n-1$, which gives the correct answers for all $n\geq 1$.
\end{remark}

     \bibliographystyle{amsplain}
     \bibliography{Boot}
\end{document}